\documentclass{svmult}
\usepackage{amssymb}
\usepackage{amsmath}
\usepackage{bbm}
\usepackage{bm}
\usepackage{booktabs}
\usepackage{multirow,subfig}

\usepackage{todonotes}

\usepackage{makeidx}         
\usepackage{graphicx}        
\usepackage{multicol}        
\usepackage[bottom]{footmisc}

\makeindex             

\begin{document}

\title*{Centering Projection Methods for Wavelet Feasibility Problems}
\author{Neil Dizon \and
Jeffrey Hogan \and Scott B. Lindstrom}
\institute{N. Dizon\\The University of Newcastle Australia\\
\email{neilkristofer.dizon@uon.edu.au}
\and \\\\J. Hogan\\The University of Newcastle Australia\\ \email{jeff.hogan@newcastle.edu.au} \and \\\\S.B. Lindstrom\\The Hong Kong Polytechnic University\\ \email{scott.b.lindstrom@polyu.edu.hk}}
%
%
\maketitle

\abstract{We revisit the feasibility approach to the construction of compactly supported smooth orthogonal wavelets on the line. We highlight its flexibility and illustrate how symmetry and cardinality properties are easily embedded in the design criteria. We solve the resulting wavelet feasibility problems using recently introduced centering methods, and we compare performance. Solutions admit real-valued compactly supported smooth orthogonal scaling functions and wavelets with near symmetry and near cardinality properties.}
\vspace{-0.3cm}

\section{Wavelet construction as a feasibility problem}
\label{sec:1}

Wavelets are traditionally constructed through \emph{multiresolution analysis} (MRA) which was introduced by Mallat \cite{mallat} and Meyer \cite{meyer}. Following MRA, Daubechies derived the first known examples of compactly supported smooth wavelets with orthonormal shifts \cite{daubechies,daubechies1}. While these wavelets have been demonstrably useful in many signal processing applications, symmetry and cardinality properties are also often desired. It is known that symmetry is incompatible with real-valuedness, orthogonality, smoothness and compact support \cite[Theorem~8.1.4]{daubechies}. In the same way, the cardinality property cannot be imposed together with all of compact support, continuity, and orthogonal shifts \cite{xia}. Recognizing these theoretical obstructions, we relax perfect symmetry or cardinality and impose only near symmetry or near cardinality. A construction technique that readily accounts for these design criteria and that easily extends to higher dimensions is preferable.

Wavelet construction has been recently formulated as a feasibility problem originally aimed at generating compactly supported smooth wavelets with orthonormal shifts \cite{fhtampaper,franklin,fhtam}. This approach handily accounts for other design criteria and allows for construction of non-tensorial wavelets in higher dimensions.
\vspace{-0.4cm}

\subsubsection*{Outline and Contributions}

In the remainder of this section, we recall one reformulation of wavelet construction as a feasibility problem. In Section~\ref{s:centeringmethods}, we recall the two centering methods we will compare: a generically proper variant of circumcentering reflections method (CRM) \cite{bauschke2018circumcentermappings,behling2018linear,behling2019convex,behling} and a new method due to Lindstrom \cite{lindstrom}. Section~\ref{s:numerical} contains our principal contribution: an experimental comparison of 2-stage global-then-local search methods, first introduced in \cite{dhlindstrom}, that combine the Douglas--Rachford method together with centering methods. This is the first such comparison for a feasibility problem, and also the first for a nonconvex problem. The results shed light on the algorithms more generally, while offering a path forward for wavelet feasibility problems specifically.
\vspace{-0.4cm}

\subsubsection*{MRA conditions and wavelet properties}

The traditional approach to the construction of wavelet orthonormal bases is based on MRA. For a more detailed discussion of the concepts that follow, refer to \cite{fhtampaper,franklin,fhtam,dhlakey,daubechies}. Henceforth, $\hat{f}$ denotes the Fourier transform of a function $f \in L^{2}(\mathbb{R},\mathbb{C})$, $\overline{A}$ is the conjugate of $A$ and denotes elementwise conjugation when $A$ is a matrix, $A[j,k]$ is the $(j,k)$-entry of a matrix $A$, and ${\rm cl}(S)$ is the closure of a set $S$.
\begin{definition}
	\label{def:MRA}
	A \textit{multiresolution analysis} for $L^2(\mathbb{R},\mathbb{C})$ consists of a sequence of closed subspaces $\{V_j\}_{j \in \mathbb{Z}}$ of $L^2(\mathbb{R},\mathbb{C})$ and a \textit{scaling function} $\varphi \in V_0$ such that the following conditions hold:
	\begin{enumerate}
		\item[(i)] the spaces $V_j$ are nested, i.e., $V_{j}\subset V_{j+1}$  for all $j \in \mathbb{Z}$, \label{def:condition1}
		\item[(ii)]$ {\rm cl}\left({\bigcup_{j \in \mathbb{Z}} V_j}\right) = L^2(\mathbb{R},\mathbb{C})$ and $\, \bigcap_{j \in \mathbb{Z}}V_j = \{0\}$, \label{def:condition2} 
		\item[(iii)] $f(\cdot) \in V_0$ if and only if $f(\cdot - k) \in V_0$ for all $k \in \mathbb{Z}$, \label{def:condition3}
		\item[(iv)] $f(\cdot) \in V_j$ if and only if $f(2(\cdot)) \in V_{j+1}$ for all $j \in \mathbb{Z}$, and \label{def:condition4}
		\item[(v)] $\{\varphi(\cdot-k)\}_{k \in \mathbb{Z}}$ forms an orthonormal basis for $V_0$. \label{def:condition5}\\
	\end{enumerate}
\end{definition}
If $\varphi$ arises from an MRA, then we are able to write $\frac{1}{2}\varphi \left(\frac{x}{2} \right) = \sum_{k \in \mathbb{Z}} h_k \varphi(x-k)$ with $\{h_k\} \in  \ell^2(\mathbb{Z})$. Taking the Fourier transforms of both sides of this equation, one obtains the \emph{scaling equation} in the Fourier domain given by $\hat{\varphi}(2\xi) = H(\xi)\hat{\varphi}(\xi)$ where $H(\xi) = \sum_{k} h_k e^{-2 \pi i k \xi}$ is the \emph{scaling filter}. Moreover, we can find a \textit{wavelet function} $\psi \in V_1 \setminus V_0$ satisfying $\frac{1}{2}\psi \left(\frac{x}{2} \right) = \sum_{k \in \mathbb{Z}} g_k \varphi(x-k)$
where $\{g_k\}_{k \in \mathbb{Z}} \in \ell^2(\mathbb{Z})$. 
Taking the Fourier transforms of both sides of this equation, one obtains $\hat{\psi}(2\xi) = G(\xi)\hat{\varphi}(\xi)$ where $G(\xi) = \sum_{k} g_k e^{-2 \pi i k \xi}$ is the \textit{wavelet filter}. If $\varphi$ has orthonormal shifts and $\left\{\psi_{j,k}:=2^{-j/2}\psi(2^{-j}x - k)\right\}_{ j,k \in \mathbb{Z}}$ forms an orthonormal basis for $L^2(\mathbb{R},\mathbb{C})$ then the wavelet matrix 
\begin{equation}U(\xi) := \begin{bmatrix}
	H(\xi) & G(\xi) \\
	H\left(\xi+\frac{1}{2}\right) & G\left(\xi +\frac{1}{2}\right) \label{def:waveletmatrix}
\end{bmatrix}
\end{equation}
is unitary for almost every $\xi \in \mathbb{R}$ and $H(0)=1$. This definition introduces a \emph{consistency condition} that $U(\xi+\frac{1}{2})=JU(\xi)$ where $J$ is the row-swap matrix.

The effectiveness of a wavelet orthonormal basis $\left\{\psi_{j,k}\right\}_{ j,k \in \mathbb{Z}}$ lies in its efficient analysis and synthesis of signals. To allow for speedy and accurate computation of the wavelet coefficients, we desire  compact support. The scaling function and wavelet are compactly supported on the interval $[0,M-1]$ if and only if  we can write $H(\xi) = \sum_{k=0}^{M-1} h_k e^{2 \pi i k \xi}$ and $G(\xi) = \sum_{k=0}^{M-1} g_k e^{2 \pi i k \xi}$ as trigonometric polynomials of degree $M-1$ \cite{franklin,fhtampaper}. Thus, we are able to write the wavelet matrix in the form 
\begin{equation*}
U(\xi)=\sum_{k=0}^{M-1}A_ke^{2\pi i k\xi} \quad \text{where}\quad A_k=\begin{bmatrix}h_k & g_k \\ (-1)^kh_k & (-1)^kgk \end{bmatrix}.
\end{equation*}

Additionally, if $\psi$ has $D$ \emph{continuous and bounded derivatives}, then this allows for better approximation using relatively fewer wavelet coefficients. Consequently, $H$, $G$ and $U$ satisfy
\begin{equation*}
\dfrac{d^kH(\xi)}{d\xi^k}\bigg|_{\xi=\frac{1}{2}}=0 \iff \dfrac{d^k G (\xi)}{d\xi^k}\bigg|_{\xi=0} = 0 \iff \left(\dfrac{d^kU(\xi)}{d\xi^k}\bigg|_{\xi=0}\right)[1,2] = 0
\end{equation*}
for all $k \in \{0,1,\ldots,D\}$, where the differentiation of the matrix is interpreted element-wise \cite{daubechies,franklin,fhtampaper}. 

Furthermore, \emph{symmetry} is another design criterion that we want $\varphi$ and $\psi$ to possess. It is known that symmetric filters applied to image processing can deal better with boundaries than asymmetric ones. A scaling function $\varphi$ is symmetric about $x=P \in (0,M-1)$ if and only if $H(\xi)=e^{4\pi iP\xi}H(-\xi)$. If $K=\mbox{diag}(-1,1) \in \mathbb{C}^{2\times 2}$, then the symmetry condition can be written in terms of the wavelet matrix as $U(\xi)=e^{4\pi i P\xi}KU(\xi)K$ \cite{dhlakey}. Note that when the scaling function is symmetric, the associated wavelet is either symmetric or anti-symmetric depending on the length of support. For conciseness, we simply say that the wavelet is symmetric.

On the other hand, \emph{cardinality} is also often sought in certain applications. A scaling function $\varphi$ is \emph{cardinal} at $P \in \mathbb{Z}$ if $\varphi(k)=\delta_{kP}$ for all $k \in \mathbb{Z}$, where $\delta$ is the Kronecker delta. A cardinal $\varphi$ admits a reconstruction formula for recovery of any function in $V_0$ from its integer samples. A necessary condition for $\varphi$ to be cardinal at $P \in \mathbb{Z}$ is $H(\xi)+(-1)^PH(\xi+\frac{1}{2})=e^{2\pi i P\xi}$ \cite{dhlakey}. Note that cardinality is desired only for the scaling function. For brevity in describing our wavelets, any mention of cardinal wavelet means that the associated scaling function is cardinal.

If we further want to guarantee that $\varphi$ and $\psi$ are \emph{real-valued}, we impose the condition that $H(\xi)=\overline{H(-\xi)}$ and $G(\xi)=\overline{G(-\xi)}$ which is equivalent to $U(\xi)=\overline{U(-\xi)}$ \cite{franklin,dhlakey}.

At this point we see that wavelet construction may be reduced to generating a matrix $U(\xi)$ satisfying the above conditions.
\vspace{-0.4cm}

\subsubsection*{Discretization by uniform sampling}

Since a trigonometric polynomial of degree $M-1$ is determined by $M$ points, we discretize $U(\xi)$ by a uniform sampling at $M$ points in $\{\frac{j}{M}\}_{j=0}^{M} \subseteq [0,1)$. By denoting each sample point by $U_j = U(\frac{j}{M})$, we form an \emph{ensemble} $\mathcal{U}:=(U_0,U_1,\dots, U_{M-1}) \in (\mathbb{C}^{2\times 2})^{M}$. The coefficient matrices $A_k$ are computed from an ensemble through an invertible $M$-point discrete Fourier transform $\mathcal{F}_M\colon (\mathbb{C}^{2\times 2})^{M} \rightarrow (\mathbb{C}^{2\times 2})^{M}:\mathcal{U} \mapsto \mathcal{A}:=(A_0,\dots,A_{M-1})$ where
\begin{align}
    A_k=({\mathcal F}_M\mathcal{U})_k&=\frac{1}{M}\sum_{j=0}^{M-1}U_je^{-2\pi ijk/M},\;\text{for}\;k\in \left\{0,1,\dots,M-1 \right\}. \label{A_k}
\end{align}

The discretized version of the consistency condition requires $U_{j +\frac{M}{2}} = J U_j$ for every $j \in \{0,1,\dots, M-1\}$. For $U(\xi)$ to be unitary almost everywhere, we need to enforce $U(\xi)$ to be unitary at $2M$ samples. Given the sample points in $\mathcal{U}$, the other set of $M$ samples may be computed to form another ensemble using $\tilde{\mathcal U}:={\mathcal F}_M^{-1}\chi_M{\mathcal F}_M({\mathcal U})$, where $(\chi_M)_j=e^{\pi ij/M}$ for $j=\{0,1,\dots, M-1\}$. Moreover, the regularity condition is imposed by forcing $(\sum_{j=0}^{M-1}j^\ell A_j)[1,2]=0$ for all $\ell \in \{0,1,\ldots,D\}$ where 
\begin{equation*}
\sum_{j=0}^{M-1}j^\ell A_j=\frac{1}{M}\sum_{k =0}^{M-1}\alpha_{\ell k}U_k 
\text{~~and~~}
\alpha_{\ell k}=\dfrac{1}{M}\sum_{j=0}^{M-1}j^\ell e^{-2\pi ik j/M}.
\end{equation*}
For the symmetry condition, we require $U_j = e^{4\pi i Pj/M}KU_{M-j}K$ for all $j \in \{1,\ldots, \frac{M}{2}\}$.  Cardinality is imposed by forcing $U_j[1,1]+(-1)^PU_{j+\frac{M}{2}}[1,1]=e^{2\pi i Pj/M}$, and the real-valuedness condition requires $U_j=\overline{U_{M-j}}$ for $j \in \{1,\ldots, \frac{M}{2}\}$. 
\vspace{-0.4cm}

\subsubsection*{The wavelet feasibility problem}

The \emph{feasibility problem} is to find a point in the intersection of a finite number of constraint sets. To reformulate wavelet construction as a feasibility problem, we treat the wavelet \emph{properties} as \emph{constraints} imposed on the discrete version of the wavelet matrix $U(\xi)$. We denote the collection of ensembles in $(\mathbb{C}^{2\times 2})^M$ that satisfy the consistency condition by $(\mathbb{C}^{2\times 2})_{J}^M$, and the collection of all $2$-by-$2$ unitary matrices by $\mathbb{U}(2)$. For an even integer $M\geq 4$ and $D=\frac{M-2}{2}$ (unless otherwise specified), we define $B_1, B_2,B_3, B_4,B_5^{(S)},B_5^{(C)} \subseteq ({\mathbb C}^{2\times 2})_J^M$ as follows.
	\begin{align*}
	B_1 &:= {\textstyle\left\{\mathcal{U}:\, U_0=\begin{pmatrix}1&0\\0&z\\ \end{pmatrix},\,|z|=1,\, U_j\in {\mathbb U}(2),\,j\in \{0,1,\dots, \frac{M}{2}\}\right\}},  \\
	B_2 &:= {\textstyle\left\{\mathcal{U}:\, ({\mathcal F}_M\chi_M({\mathcal F}_M)^{-1}(\mathcal{U}))_j\in {\mathbb U}(2),\, j\in \{0,1,\dots, \frac{M}{2}\}\right\}},  \\
	B_3 &:= {\textstyle\left\{\mathcal{U}:\, \left(\sum_{j=0}^{M-1}\alpha _{\ell k}U_k\right)[1,2]=0,\ 0\leq\ell\leq D\right\}},\\
	B_4 &:= \textstyle{\left\{\mathcal{U}: U_j = \overline{U_{M-j}}, \, j \in \{1,2,\dots, \frac{M}{2}\} \right\}}, \\
	B_5^{(S)}&:= \left\{\mathcal{U}: \|U_j-e^{2\pi i Pj/M}KU_{M-j}K\|<\gamma, \, j \in \{1,2,\dots, M/2\} \right\}\\
	B_5^{(C)}&:={\textstyle \left\{\mathcal{U}:\left\|U_j[1,1]+(-1)^PU_{j+\frac{M}{2}}[1,1]-e^{2\pi i Pj/M}\right\|<\gamma, \, j \in \{1,2,\dots, \frac{M}{2}\} \right\}}.
	\end{align*}

Note that $B_1$ and $B_2$ are nonconvex constraint sets that correspond to the unitarity condition at $2M$ sample points. The subspaces $B_3$ and $B_4$ are constraint sets for regularity and real-valuedness, respectively. Moreover, $B_5^{(S)}$ and $B_5^{(C)}$ are convex sets that promote near symmetry and near cardinality properties, respectively. Notice the introduction of a small positive number $\gamma$ in the definition of $B_5^{(S)}$ and $B_5^{(C)}$ to get around the theoretical obstructions for obtaining perfect symmetry and cardinality \cite{dhlakey}. In summary, we have the following feasibility problems.

\begin{problem}[Nearly symmetric wavelets]\label{prob:waveletproblem1}
	The feasibility problem for constructing compactly supported real-valued smooth nearly symmetric orthogonal wavelets is to find an ensemble $\ \mathcal{U}\in \bigcap_{k=1}^4 B_k \cap B_5^{(S)} \subseteq ({\mathbb C}^{2\times 2})_J^M$.
\end{problem}

\begin{problem}[Nearly cardinal wavelets]\label{prob:waveletproblem2}
	The feasibility problem for constructing compactly supported real-valued smooth nearly cardinal orthogonal wavelets  is to find an ensemble $\ \mathcal{U}\in \bigcap_{k=1}^4 B_k \cap B_5^{(C)} \subseteq ({\mathbb C}^{2\times 2})_J^M$.
\end{problem}
\vspace{-0.4cm}

\section{Centering  methods for feasibility problems}\label{s:centeringmethods}

The original works that solved wavelet feasibility problems for compactly supported smooth orthogonal wavelets employed the \emph{Douglas--Rachford (DR) algorithm} \cite{drachford,LSsurvey} to solve \emph{Pierra's product space reformulation} \cite{pierra} of the feasibility problem. The method demonstrated surprising robustness in this context  \cite{fhtampaper,franklin,fhtam,ddhtam}. Convergence plots frequently feature the tell-tale characteristics of local spiraling during convergence; such features are described in \cite{lindstrom}. The spiraling is associated with longer runs for numerical implementations \cite{franklin} and presents an opportunity to accelerate convergence \cite{lindstrom}.

In this section, we recall the DR operator, the generalized circumcentered reflections method operator (GCRM) \cite{dhlindstrom,behling} and the new centering operator $L_T$ introduced by Lindstrom in \cite{lindstrom}. We expect the two centering methods to accelerate convergence to feasible solutions.

For a closed subset $C$ of a Hilbert space $\mathcal{H}$, we define the operator $P_C:\mathcal{H}\to C$ by $P_C x \in {\rm argmin}_{z \in C}\|z-x\|$; it is a selector for the closest point projection for $C$. Its associated \emph{reflector} is defined as $R_V:=2P_V-Id$ where $Id$ is the identity map. Given three points $x,y,z \in \mathcal{H}$, we denote $C(x,y,z)$ to be their \emph{circumcenter}, which is equidistant to the given points and lies on the affine subspace they define. The circumcenter exists whenever $x,y,z$ are not simultaneously distinct and colinear; for more on existence and formulae for computation, see \cite{bauschke2018circumcentermappings,bauschke2018circumcenters}.

\begin{definition}
Let  $V$ and $W$ be nonempty subsets of $\mathcal{H}$.
\begin{enumerate}
	\item The \emph{DR operator} for $V$ and $W$ is defined as $T(x):=x-P_V(x) +P_WR_V(x)$.
	\item The \emph{circumcentering reflections method} operator is defined as $ CRM(x):=C(x,R_V(x),R_WR_V(x))$. For history and properties, see \cite{behling2018linear,behling2019convex,behling}. 
	\item The \emph{GCRM operator} is defined as
	\begin{equation*}
	C_{V,W}(x) := \begin{cases}
	T(x) & \text{if}\; x,R_Vx,R_WR_Vx\;\;\text{are colinear;}\\
	CRM(x) & \text{otherwise.}
	\end{cases}
	\end{equation*}
	\item The centering operator $L_T$ from \cite{lindstrom} is defined as: 
		\begin{equation*}
		L_{T}(x) := \begin{cases}
		C(x,2Tx-x,\pi_T(x)) & \text{if}\; x,2Tx-x,\pi_T(x)\;\;\text{are not colinear;}\\
		T^2x & \text{otherwise.}
		\end{cases}
		\end{equation*}
		where $\pi_T(x)=2(T^2x-Tx)+2P_{\text{span}(T^2x-Tx)}(Tx-x)+x$.
\end{enumerate} 
\end{definition}

Lindstrom discovered that for some prototypical feasibility problems for which Lyapunov functions are known, CRM returns the minimizer of a quadratic surrogate for the local Lyapunov function \cite{lindstrom}. Lindstrom showed that $L_T$'s lack of dependence on subproblems (in our case, reflections) allows it to recapture this property in settings where CRM may not, such as for the primal-dual implementation of ADMM/Douglas--Rachford for basis pursuit. In our setting, this possible improvement in stability carries the computational cost that one application of $L_T$ requires two applications of the pair of projections $P_V$ and $P_W$, instead of just one pair for CRM.

For numerical implementations, we set up a 2-stage DR-GCRM and a 2-stage DR-$L_T$. In stage 1, we exploit the greater global robustness of DR to find local basins of attraction to feasible points, and thereafter, in stage 2, we apply centering methods to obviate local spiraling thereto. It has already been shown experimentally that this approach consistently outperforms a full run of DR in the context of solving wavelet feasibility problems \cite{dhlindstrom}. In the next section, we use GCRM  and $L_T$ as the local methods of 2-stage global-then-local search algorithms, in order to solve  Problem \ref{prob:waveletproblem1} and \ref{prob:waveletproblem2}.
\vspace{-0.4cm}

\section{Numerical Results}\label{s:numerical}

We use a product space technique similar to those employed in \cite{fhtam,fhtampaper,franklin,dhlindstrom} to convert our many-set feasibility problems into 2-set problems amenable to solution by the methods described above.\\

\noindent \textit{Problem~\ref{prob:waveletproblem1}}: The constraints for obtaining nearly symmetric wavelets are
\begin{align*}
V &:=B_1\times B_2 \times \left(B_3\cap B_4\right) \times B_5^{(S)} \subseteq \left((\mathbb{C}^{2\times 2})_{J}^M\right)^4,\\
W &:= \left\{(\mathcal{U}_j)_{j=1}^{4} \in \left((\mathbb{C}^{2\times 2})_{J}^M\right)^4: \mathcal{U}_1 = \mathcal{U}_2 = \mathcal{U}_3 = \mathcal{U}_4\right\}.
\end{align*}

\noindent \textit{Problem~\ref{prob:waveletproblem2}}: The constraints for obtaining nearly cardinal wavelets are
\begin{align*}
V &:=B_1\times B_2 \times \left(B_3\cap B_4\right) \times B_5^{(C)} \subseteq \left((\mathbb{C}^{2\times 2})_{J}^M\right)^4,\\
W &:= \left\{(\mathcal{U}_j)_{j=1}^{4} \in \left((\mathbb{C}^{2\times 2})_{J}^M\right)^4: \mathcal{U}_1 = \mathcal{U}_2 = \mathcal{U}_3 = \mathcal{U}_4\right\}.
\end{align*}

The projection of a $4$-tuple of ensembles onto $W$ is obtained by averaging the $4$ ensembles. Notice that the set $V$ and its projection $P_V$ are different for the two problems, though this should create no confusion because we will only discuss one problem at a time. We have $P_V=P_{B_1}\times P_{B_2} \times (P_{B_3}P_{B_4}) \times P_{B_5^{(\eta)}}$, where $\eta$ is $C$ or $S$ respectively for the two problems. Because $P_{B_3}(B_4) \subset B_4$ and $B_4$ is a subspace, the identity  $P_{B_3}P_{B_4}= P_{B_3 \cap B_4}$ admits the \emph{constraint-reduction reformulation} we have used; see \cite{ddhtam}.

In what follows, we solve Problem \ref{prob:waveletproblem1} and \ref{prob:waveletproblem2} with  $M=6, D=1$ and $\gamma=0.5$. We compare the performance of 2-stage DR-GCRM with 2-stage DR-$L_T$. We initialize at $100$ random ensembles that satisfy the consistency condition. Throughout, we let $(x_n)_{n\in \mathbb{N}}$ be the sequence of iterates generated by the projection algorithm under consideration. We fix a tolerance $\varepsilon :=10^{-9}$ and use the stopping criterion $\varepsilon_n:=\|P_VP_W(x_n) - P_W(x_n)\| < \varepsilon$, whereupon $P_W(x_n)$ is a feasible point. In implementing a 2-stage method, we first run DR until the gap distance $\varepsilon_n$ reaches a $10^{-2}$ threshold; thereafter we switch to applying GCRM or $L_T$. We declare a particular run to have \emph{solved} the feasibility problem whenever it attains the threshold of $\varepsilon$ within $20,000$ iterations. We provide statistics on the \emph{number of iterations} needed, which is our main performance measure. We do not report the number of iterates required for $\varepsilon_n$ to obtain the threshold $10^{-2}$, because it is the same for both 2-stage algorithms. We only report the number of iterates needed thereafter.

\begin{figure}[h]
	\centering
	\subfloat{\includegraphics[width = 0.4\linewidth]{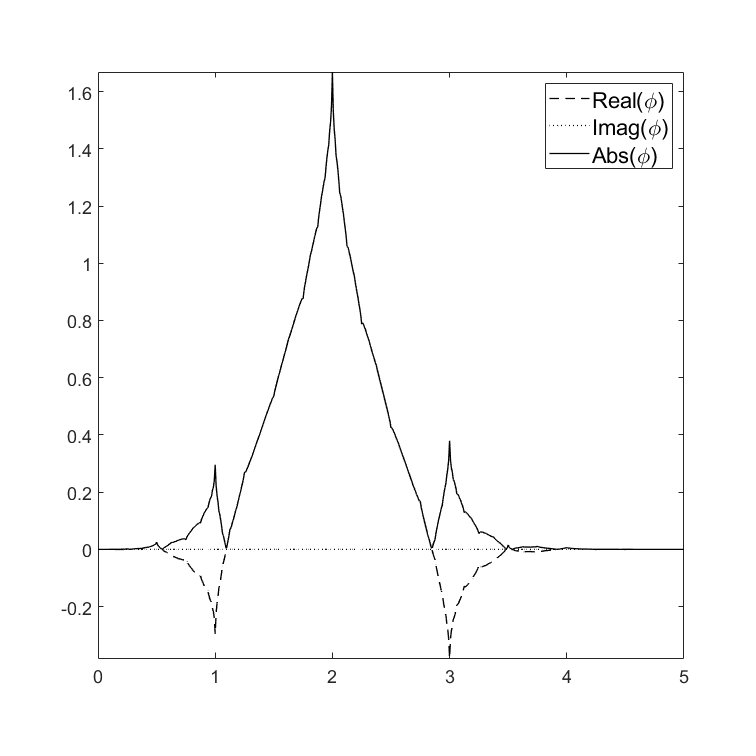}}
	\subfloat{\includegraphics[width = 0.4\linewidth]{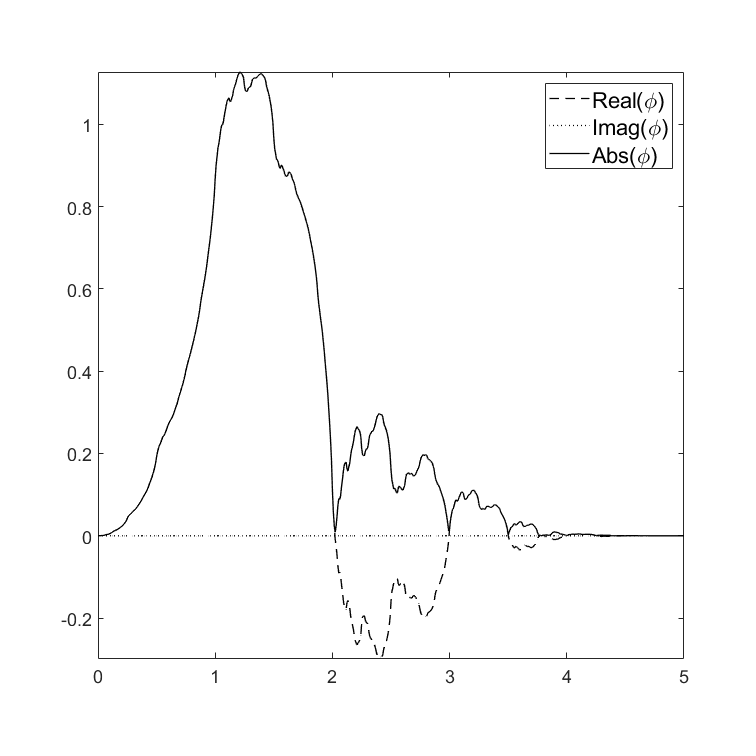}}
	\caption{Nearly symmetric (about $x=2$) and nearly cardinal (at $x=1$) scaling functions plotted by employing \emph{cascade algorithm} on filters generated from the feasible ensembles solved in Problem \ref{prob:waveletproblem1} and \ref{prob:waveletproblem2}, respectively.}  \label{fig:symmetricexample}
\end{figure}

\begin{table}[h]
	\centering
	\begin{tabular}{@{}clccccccc@{}}
		\toprule
		\multirow{2}{*}{} & \multicolumn{1}{c}{\multirow{2}{*}{Algorithm}} & \multirow{2}{*}{\begin{tabular}[c]{@{}c@{}}cases\\ solved\end{tabular}} & \multicolumn{1}{l}{\multirow{2}{*}{\begin{tabular}[c]{@{}l@{}}solved\\ by all\end{tabular}}} & \multicolumn{5}{c}{when solved by all} \\ \cmidrule(l){5-9} 
		& \multicolumn{1}{c}{} &  & \multicolumn{1}{l}{} & \multicolumn{1}{l}{wins} & \multicolumn{1}{l}{Q1} & \multicolumn{1}{l}{mean} & \multicolumn{1}{l}{Q3} & \multicolumn{1}{l}{median} \\ \midrule
		\multirow{3}{*}{Problem 1} & DR & 51 & 51 & 0 & 194 & 211 & 215 & 201 \\
		& GCRM & 51 & 51 & 13 & 36 & 28 & 40 & 38 \\
		& $L_T$ & 51 & 51 & 38 & 29 & 36 & 39 & 33 \\ \midrule
		\multirow{3}{*}{Problem 2} & DR & 96 & 79 & 0 & 176 & 182 & 186 & 185 \\
		& GCRM & 79 & 79 & 22 & 31 & 33 & 35 & 33 \\
		& $L_T$ & 96 & 79 & 57 & 28 & 32 & 33 & 31 \\ \bottomrule
	\end{tabular}
	\caption{Performance during stage 2 of a 2-stage search.}
	\label{table:probs}
\end{table}

Table~\ref{table:probs} summarizes the numerical results. $L_T$ solved every problem DR solved. For Problem~\ref{prob:waveletproblem2}, GCRM was less stable than $L_T$, which is consistent with what one might expect, given that $L_T$ is constructed to retain the property of minimizing a surrogate Lyapunov function in situations where GCRM's dependence on subproblems may cause instability \cite{lindstrom}. Interestingly, for Problem~\ref{prob:waveletproblem1}, GCRM also solved every problem DR solved. When both algorithms converged, $L_T$ and GCRM performed quite similarly, which is what one would expect if both methods are constructing, from their respective sampling points, relatively similar quadratic surrogates for the underlying Lyapunov function. However, one should remember that computing a single centering step of $L_T$ requires computing twice the number of projection substeps that are needed by a single step of GCRM.
\vspace{-0.3cm}

\section{Conclusion}
We have  shown how the symmetry and cardinality constraints are readily accounted for in the feasibility approach to wavelet construction. Numerical results also shed light on local behaviour of $L_T$ and GCRM. We speculate that both are viable heuristics that may be applied to deal with wavelet feasibility problems for higher dimensional constructions,  and we suggest this as the next step of research.

\bibliographystyle{siam}
\bibliography{bib}
%



\printindex
\end{document}